\documentclass{article}
\usepackage{amssymb,latexsym,amsmath,amsthm,mathrsfs}
\usepackage{graphicx}
\newcommand{\comment}[1]{}

\begin{document}
\title{On amicable numbers\footnote{Originally published as
{\em De numeris amicabilibus},
Nova acta eruditorum (1747), 267--269.
E100 in the Enestr{\"o}m index.
Translated from the Latin by Jordan Bell,
Department of Mathematics, University of Toronto, Toronto, Ontario, Canada.
Email: jordan.bell@gmail.com}}
\author{Leonhard Euler}
\date{}
\maketitle

At this time in which mathematical Analysis has opened the way
to many profound observations, those problems which have to do with the
nature and properties of numbers seem almost completely neglected
by Geometers, and the contemplation of numbers has been judged by many
to add nothing to Analysis. Yet truly the investigation of the properties of
numbers on many occasions requires more acuity than the subtlest questions
of geometry, and for this reason it seems improper to
neglect arithmetic questions for those. And indeed the greatest
thinkers who are recognized as having made the most important contributions
to Analysis have judged the affection of numbers as not unworthy,
and in pursuing them have expended much work and study.
Namely, it is known that Descartes, even though occupied with the most
important meditations on both universal Philosophy and especially
Mathematics, spent no little effort uncovering amicable numbers;
this matter was then pursued even more by van Schooten.\footnote{Translator: The {\em Opera omnia} cites many relevant works. Also, the {\em Opera omnia}
refers to Euler's papers E152 and E798, both on amicable numbers.}
One calls two numbers amicable if the sum of the aliquot parts of each
yields the other;
220 and 284 are numbers of this type; for the aliquot parts of the first
220, that is the divisors less than the number itself, 
\[
1+2+4+5+10+11+20+22+44+55+110,
\]
produce the sum 284, and the aliquot parts of this number 284,
\[
1+2+4+71+142,
\]
produce in turn 220. There is also no doubt that beside these two numbers
many others and even infinitely many may be given which have the same
property; but neither Descartes nor after him van Schooten
were to uncover more than these three pairs of numbers,
 even though it
seems that they expended not a small amount of study for finding many.
And indeed the method which both of them used is such that
they could scarcely find more pairs of amicable numbers;\footnote{Translator: cf. {\em Notes on Th\=abit ibn Qurra and his rule for
amicable numbers}, Historia Math. \textbf{16} (1989), no. 4, 373--378.}
for they took the numbers to be contained in the forms
$2^nxy$ and $2^nz$, where $x,y$ and $z$ denote prime numbers that
should be such that first $z=xy+x+y$, and then
too $2^n(x+y+2)=xy+x+y+1$.\footnote{Translator: Namely, if $x,y,z,n$ satisfy
these conditions then
$2^nxy$ and $2^nz$ are an amicable pair.}
Thus they successively took different values for the exponent $n$ and
for each sought prime numbers $x$ and $y$ so that the previous equation
is satisfied; as soon as it was the case that $xy+x+y$ yielded a prime
number, the given formulas $2^nxy$ and $2^nz$ exhibited amicable numbers.
And it is easy to see that by proceeding this way to larger exponents $n$, soon
such large numbers $xy+x+y$ are reached that one can no longer discern
whether they are prime or not, since the table of prime numbers
has not yet been extended beyond 100000.

It is clear that this question is restricted no small amount beyond necessity
as long as amicable numbers are assumed to be included only in the 
formulas given so far. 
I had contemplated this carefully, and by calling on the aid of several
artifices taken from the nature of divisors I obtained 
many other pairs of amicable numbers, thirty of which including the three
already known I shall communicate here; also, so that their origin and nature
can be clearly perceived, I shall write them in factored form.
Thus the amicable numbers are:\footnote{Translator: Euler does not explain how he found these amicable pairs. Pair XIII is not in fact an amicable pair.
See Ed Sandifer's November 2005 {\em How Euler did it}, www.maa.org/news/howeulerdidit.html}

\[
\begin{array}{rlrl}
\textrm{I.}&\begin{cases}2^2\cdot 5\cdot 11\\2^2\cdot 71\end{cases}
&\textrm{XVI.}&\begin{cases}2^5\cdot 37\cdot 12671\\ 2^5\cdot 227\cdot 2111\end{cases}\\
\textrm{II.}&\begin{cases}2^4\cdot 23\cdot 47\\2^4\cdot 1151\end{cases}
&\textrm{XVII.}&\begin{cases}2^5\cdot 53\cdot 10559\\2^5\cdot 79\cdot 7127\end{cases}\\
\textrm{III.}&\begin{cases}2^7\cdot 191\cdot 383\\2^7\cdot 73727\end{cases}
&\textrm{XVIII.}&\begin{cases}2^6\cdot 79\cdot 11087\\2^6\cdot 383\cdot 2309\end{cases}\\
\textrm{IV.}&\begin{cases}2^2\cdot 23\cdot 5\cdot 137\\2^2\cdot 23\cdot 827\end{cases}
&\textrm{XIX.}&\begin{cases}2^2\cdot 11\cdot 17\cdot 263\\2^2\cdot 11\cdot 43\cdot 107\end{cases}\\
\textrm{V.}&\begin{cases}3^2\cdot 5\cdot 13\cdot 11\cdot 19\\3^2\cdot 5\cdot 13\cdot 239\end{cases}
&\textrm{XX.}&\begin{cases}3^3\cdot 5\cdot 7\cdot 71\\3^3\cdot 5\cdot 17\cdot 31\end{cases}\\
\textrm{VI.}&\begin{cases}3^2\cdot 7\cdot 13\cdot 5\cdot 17\\3^2\cdot 7\cdot 13\cdot 107\end{cases}
&\textrm{XXI.}&\begin{cases}3^2\cdot 5\cdot 13\cdot 29\cdot 79\\3^2\cdot 5\cdot 13\cdot 11\cdot 199\end{cases}\\
\textrm{VII.}&\begin{cases}3^2\cdot 7^2\cdot 13\cdot 5\cdot 41\\3^2\cdot 7^2\cdot 13\cdot 251\end{cases}
&\textrm{XXII.}&\begin{cases}3^2\cdot 5\cdot 13\cdot 19\cdot 47\\3^2\cdot 5\cdot 13\cdot 29\cdot 31\end{cases}\\
\textrm{VIII.}&\begin{cases}2^2\cdot 5\cdot 131\\2^2\cdot 17\cdot 43\end{cases}
&\textrm{XXIII.}&\begin{cases}3^2\cdot 5\cdot 13\cdot 19\cdot 37\cdot 1583\\
3^2\cdot 5\cdot 13\cdot 19\cdot 227\cdot 263\end{cases}\\
\textrm{IX.}&\begin{cases}2^2\cdot 5\cdot 251\\2^2\cdot 13\cdot 107\end{cases}
&\textrm{XXIV.}&\begin{cases}3^3\cdot 5\cdot 31\cdot 89\\3^3\cdot 5\cdot 7\cdot 11\cdot 29\end{cases}\\
\textrm{X.}&\begin{cases}2^3\cdot 17\cdot 79\\2^3\cdot 23\cdot 59\end{cases}
&\textrm{XXV.}&\begin{cases}2\cdot 5\cdot 7\cdot 60659\\2\cdot 5\cdot 23\cdot 29\cdot 673\end{cases}\\
\textrm{XI.}&\begin{cases}2^4\cdot 23\cdot 1367\\2^4\cdot 53\cdot 607\end{cases}
&\textrm{XXVI.}&\begin{cases}2^3\cdot 31\cdot 11807\\2^3\cdot 11\cdot 163\cdot 191\end{cases}\\
\textrm{XII.}&\begin{cases}2^4\cdot 17\cdot 10303\\2^4\cdot 167\cdot 1103\end{cases}
&\textrm{XXVII.}&\begin{cases}3^2\cdot 7\cdot 13\cdot 23\cdot 79\cdot 1103\\3^2\cdot 7\cdot 13\cdot 23\cdot 11\cdot 19\cdot 367\end{cases}\\
\textrm{XIII.}&\begin{cases}2^4\cdot 19\cdot 8563\\2^4\cdot 83\cdot 2039\end{cases}
&\textrm{XVIII.}&\begin{cases}2^3\cdot 47\cdot 2609\\2^3\cdot 11\cdot 59\cdot 173\end{cases}\\
\textrm{XIV.}&\begin{cases}2^4\cdot 17\cdot 5119\\2^4\cdot 239\cdot 383\end{cases}
&\textrm{XXIX.}&\begin{cases}3^3\cdot 5\cdot 23\cdot 79\cdot 1103\\3^3\cdot 5\cdot 23\cdot 11\cdot 19\cdot 367\end{cases}\\
\textrm{XV.}&\begin{cases}2^5\cdot 59\cdot 1103\\2^5\cdot 79\cdot 827\end{cases}
&\textrm{XXX.}&\begin{cases}3^2\cdot 5^2\cdot 11\cdot 59\cdot 179\\3^2\cdot 5^2\cdot 17\cdot 19\cdot 359\end{cases}
\end{array}
\]

\end{document}